%% file: main.tex
\documentclass[12pt,leqno]{article}

\usepackage[T2A]{fontenc}
\usepackage{amsfonts}
\usepackage{amsmath}
\usepackage{amssymb}
\usepackage{comment}
\usepackage[dvips]{graphicx}
\usepackage{subfigure}

\newtheorem{theorem}{Theorem}[section]
\newtheorem{lemma}[theorem]{Lemma}

\newtheorem{remark}[theorem]{Remark}

\newcommand{\set}[1]{\left \{ #1 \right \}}                     
\newcommand{\setst}[2]{\left \{ #1 \mid #2 \right \}}           
\newcommand{\abs}[1]{\left| #1 \right|}

\providecommand{\R}{\mathbb{R}}

\renewcommand{\phi}{\varphi}

\renewcommand{\bar}{\overline}

\newcommand{\calB}{\mathcal{B}}
\newcommand{\calP}{\mathcal{P}}
\newcommand{\calF}{\mathcal{F}}

\newcommand{\deltain}{\delta^{\rm in}}
\newcommand{\deltaout}{\delta^{\rm out}}

\newenvironment{proof}{\par\noindent%
{\bf Proof.\par\nopagebreak}}{\unskip\nobreak\enskip$\square$\par\bigskip}

\newenvironment{numitem}{%
   \refstepcounter{equation}%
   \begin{enumerate}%
      \item[(\theequation)]
}{%
   \end{enumerate}
}

\newcommand{\refeq}[1]{(\ref{eq:#1})}                
\newcommand{\reffig}[1]{Fig.~\ref{fig:#1}}           
\newcommand{\refth}[1]{Theorem~\ref{th:#1}}          
\newcommand{\reflm}[1]{Lemma~\ref{lm:#1}}            
\newcommand{\refsec}[1]{Section~\ref{sec:#1}}        
\newcommand{\refrem}[1]{Remark~\ref{rem:#1}}         

\makeindex

\begin{document}

\title
{
   On Ear Decompositions of \\
   Strongly Connected Bidirected Graphs
}

\author
{
   Maxim A. Babenko
   \thanks
   {
      Dept. of Mechanics and Mathematics,
      Moscow State University, Vorob'yovy Gory, 119899 Moscow,
      Russia, \textsl{email}: mab@shade.msu.ru.
      Supported by RFBR grants 03-01-00475, 05-01-02803, and 06-01-00122.
   }
}

\maketitle

\begin{abstract}
   \input abstract.tex
\end{abstract}



\section{Introduction}
\label{sec:intro}

For an arbitrary undirected graph~$G$ we write $V_G$
(resp. $E_G$) to denote the set of nodes (resp. edges) of~$G$.
In case $G$ is directed we speak of arcs rather than edges and write~$A_G$
instead of~$E_G$. The same notation will be used for walks, paths, cycles etc.

Consider a digraph $G$ and its arbitrary subgraph $H$ (that is, $V_H \subseteq V_G$,
$A_H \subseteq A_G$). An \emph{ear} of $H$ w.r.t.~$G$ is a path~$P$ in~$G$
such that: (i) both ends of~$P$ are in~$V_H$; (ii) no inner node of~$P$ is in~$V_H$;
(iii) $A_P \cap A_H = \emptyset$.
In particular, an ear can consist of a single arc~$a$ with both head and tail nodes in~$V_H$;
as long as this is not confusing we denote this ear by~$a$.
By $H' := H + P$ we denote a new digraph with $V_{H'} := V_H \cup V_P$, $A_{H'} := A_H \cup A_P$.
Also, for a collection~$\calP$ of ears we denote by $G + \calP$ the result of adding all ears
from~$\calP$ to~$G$.

Recall that a digraph~$G$ is called \emph{strongly connected} if for any pair of nodes in~$G$
the former one is reachable from the latter by a path or, equivalently,
the underlying undirected graph of~$G$ is connected and each arc of~$G$ is contained in
a cycle.

For a pair of strongly connected digraphs $G, H$, where $H$ is a subgraph of~$G$,
we define an \emph{ear decomposition} of $G$ starting from~$H$ to be
a sequence of strongly connected subgraphs of~$G$
$$
   H = G_0, G_1, \ldots, G_{k-1}, G_k = G,
$$
where $G_{i+1}$ is obtained from $G_i$ by adding an ear of~$G_i$ w.r.t.~$G$
($0 \le i < k$). Clearly, an ear decomposition is not unique.

\begin{remark}
\label{rem:dir_ear_decomp}
   One can easily see that the requirement for $G_1, \ldots, G_k$ to be
   strongly connected can be dropped since adding an ear to a strongly connected
   digraph preserves strong connectivity.
   This will not be the case for the class of bidirected graphs so we keep this
   requirement to make our definitions more symmetric.
\end{remark}

A central fact about ear decompositions of strongly connected digraphs
is stated in the next folklore theorem:
\begin{theorem}
\label{th:dir_decomp}
   For any strongly connected digraph~$G$ and an arbitrary strongly connected
   subgraph~$H$ of $G$ there exists an ear decomposition of~$G$ starting from~$H$.
\end{theorem}

\medskip

The main goal of this paper is to extend the notion of ear decomposition
and \refth{dir_decomp} to the class of \emph{bidirected} graphs.
It turns out that this generalization will naturally contain certain
well-known decomposition results from matching theory.

The notion of bidirected graphs was introduced by
Edmonds and Johnson~\cite{EJ-70} in connection with one important class of integer linear programs
generalizing problems on flows and matchings; for
a survey, see also~\cite{law-76,sch-03}.

Recall that in a \emph{bidirected} graph $G$ three types of edges are allowed:
(i)~a~standard directed edge, or an \emph{arc}, that leaves one node and enters another one;
(ii)~a~nonstandard edge leaving both of its ends; or
(iii)~a~nonstandard edge entering both of its ends.

When both ends of an edge coincide, the edge becomes a \emph{loop}.

We borrow the notation that was introduced for undirected graphs and
write~$V_G$ (resp. $E_G$) to denote the set of nodes (resp. edges) of
a bidirected graph~$G$.

A \emph{walk} in a bidirected graph $G$ is an alternating sequence
$P = (s = v_0, e_1, v_1, \ldots, e_k, v_k = t)$ of nodes and edges such that
each edge $e_i$ connects nodes $v_{i-1}$ and $v_i$, and for
$i = 1, \ldots, k-1$, the edges $e_i,e_{i+1}$ form a \emph{transit pair}
at $v_i$, which means that one of $e_i,e_{i+1}$ enters and the
other leaves~$v_i$. Note that $e_1$ may enter $s$ and $e_k$ may
leave $t$; nevertheless, we refer to $P$ as a walk from $s$ to $t$,
or an $s$--$t$ \emph{walk}. $P$ is \emph{cyclic} if $v_0=v_k$ and
the pair $e_1,e_k$ is transit at $v_0$; cyclic walks are usually
considered up to cyclic shifts. Observe that an $s$--$s$ walk
is not necessarily cyclic.

A walk is called \emph{edge-simple} (or a \emph{path}) if all its edges are different.
If $v_i \ne v_j$ for all $1 \le i < j < k$ and $1 < i < j \le k$,
then walk~$P$ is called \emph{node-simple} (or a \emph{simple} path).
Note that the ends of a simple path need not be distinct.
As usually, a cyclic edge-simple walk is called a \emph{cycle}. A node-simple cyclic
walk is called a \emph{simple} cycle.

We now extend the notions of strong connectivity and ear decomposition to the
class of bidirected graphs. We call a bidirected graph~$G$ \emph{strongly connected}
if its underlying undirected graph is connected and each edge of $G$ is
contained in a cycle.

For a bidirected graph~$G$ and its subgraph~$H$ an \emph{ear} of $H$ w.r.t.~$G$ is
a path~$P$ in~$G$ such that: (i) both ends of~$P$ are in~$V_H$;
(ii) no inner node of~$P$ is in~$V_H$;
(iii) $E_P \cap E_H = \emptyset$.
As earlier, we use notation~$e$ to denote the ear consisting of a single edge~$e$.

\begin{figure}[tb]
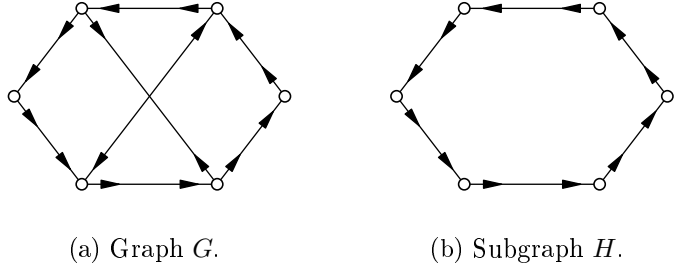

    \centering
    \subfigure[Graph~$G$.]{
      \includegraphics{pics/twoedges.2}%
    }
    \hspace{1cm}%
    \subfigure[Subgraph~$H$.]{
      \includegraphics{pics/twoedges.1}%
    }
    \caption{A pair of strongly connected bidirected graphs $G$, $H$ such that
    $H$ is a subgraph of $G$ and $G$ cannot be obtained from $H$ by adding a single ear.}
    \label{fig:twoedges}
\end{figure}
One can see that unlike the case of directed graphs adding an ear to a strongly connected
instance may produce a graph that is not strongly connected (cf.~\refrem{dir_ear_decomp}).
Moreover, being restated in terms of bidirected graphs, \refth{dir_decomp} becomes false.
To see this, consider an example depicted in \reffig{twoedges}. Both graphs $G$, $H$ are
strongly connected and $H$ can be obtained from $G$ by adding two edges.
However, adding only one of these edges does not produce a strongly connected instance.

To overcome this obstacle one needs to allow a pair of ears to be added on certain
steps. More formally, consider strongly connected bidirected graphs $G$ and $H$
such that $H$ is a subgraph of $G$. Also, consider a collection $P_1, \ldots, P_k$
of ears of~$H$ w.r.t.~$G$. We denote by $H' := H + P_1 + \ldots + P_k$ the result
of adding all ears~$P_i$ to $H$. In particular,
$$
   V_{H'} := V_H \cup V_{P_1} \cup \ldots \cup V_{P_k}, \qquad
   E_{H'} := E_H \cup E_{P_1} \cup \ldots \cup E_{P_k}.
$$
Consider a strongly connected bidirected graph $G$ and its strongly connected
subgraph~$H$. An \emph{ear decomposition} of $G$ starting from~$H$ is
a sequence of strongly connected subgraphs of~$G$
$$
   H = G_0, G_1, \ldots, G_{k-1}, G_k = G,
$$
where $G_{i+1}$ is obtained from $G_i$ by adding a single ear of~$G_i$ w.r.t.~$G$
or an edge-disjoint pair of such ears ($0 \le i < k$). In case $G_{i+1}$ is obtained
from $G_i$ by adding only one ear we call it a \emph{single-ear step}; otherwise
we are referring to it as a \emph{double-ear step}.

The required generalization of~\refth{dir_decomp} can now be stated as follows:
\begin{theorem}
\label{th:bidir_decomp}
   For any strongly connected bidirected graph~$G$ and an arbitrary strongly connected
   subgraph~$H$ of $G$ there exists an ear decomposition of~$G$ starting from~$H$.
\end{theorem}

The rest of the paper is organized as follows. In~\refsec{two_edges} we prove
a certain special case of~\refth{bidir_decomp} (that may be interesting for its own
sake). Sections~\ref{sec:ss} and \ref{sec:barrier} contain
some basic results regarding the so-called skew-symmetric graphs,
which are used later in \refsec{bidir_decomp_proof}, where a complete proof of~\refth{bidir_decomp}
is given. In~\refsec{applications} we show how Two Ears Theorem on matching covered
graphs can be derived from our results.

\section{Two Edges Theorem}
\label{sec:two_edges}

\begin{theorem}
\label{th:two_edges}
   Let $G$ be a strongly connected bidirected graph with all edges standard;
   let~$E$ be a nonempty collection of bidirected edges with both ends in $V_G$
   such that each edge in $E$ is nonstandard and $G + E$ is strongly connected.
   Then there exist a pair of edges $e_1, e_2 \in E$
   such that $G + e_1 + e_2$ is also strongly connected.
\end{theorem}

Suppose towards contradiction that there exists a graph~$G$
and a collection of nonstandard edges~$E$ such that $\abs{E} > 2$ and $G + E$ is strongly-connected
but no proper subset $E' \subset E$ satisfies this property.
In what follows we regard $G$ as a standard directed graph denoting the set of its arcs by~$A_G$.
Each edge in $E$ is of two possible kinds: it either enters both ends or
leaves them; according to this, we divide $E$ into the subsets $E^+$ and $E^-$
respectively.

Consider a cycle~$C$ in $G + E$ that uses at least one nonstandard edge.
Then, $C$ traverses equal number of edges from $E^+$
and $E^-$. By assumption of minimality of $E$, $C$ traverses
all edges of $E^+$ and $E^-$, and hence $\abs{E^+} = \abs{E^-}$.
Put
$$
   E^+ = \set{e^+_1, \ldots, e^+_n}, \quad
   E^- = \set{e^-_1, \ldots, e^-_n}.
$$

We transform $G$ and $E$ in order to make
sure that all ends of edges in~$E$ are distinct.
To this aim we do the following:
(i) split each node~$v \in V_G$ into a sufficient number of pairs
$v_i^+, v_i^-$;
(ii) for each node $v \in V_G$ add arcs $(v_i^+, v_j^-)$ between all possible pairs;
(iii) transform each arc $(u, v) \in A_G$ into a collection of arcs $(u^-_i, v^+_j)$
going between all possible pairs.
Clearly, this transformation preserves strong connectivity of~$G$.

Finally, each edge $\set{u,v} \in E^+$ (resp.~$\set{u,v} \in E^-$)
is transformed into an edge $\set{u^+_i, v^+_j}$
(resp.~$\set{u^-_i, v^-_j}$) of the same type. Here we choose
``fresh'' values of $i, j$ for each edge to guarantee
that all ends are distinct. In what follows
we keep the same notation $G$ and $E$ to denote the resulting
graph and the resulting set of nonstandard edges.

\medskip

Recall \cite{sch-03} that for a given nonempty set~$V$ a pair $(X,Y)$,
$X, Y \subseteq V$, is said to be~\emph{crossing} if
$X \cap Y \ne \emptyset$, $X \cup Y \ne V$,
$X \setminus Y \ne \emptyset$, and $Y \setminus X \ne \emptyset$.
A~family of sets $\calF \subseteq 2^V$ is called
\emph{crossing} if $X \cap Y, X \cup Y \in \calF$
for every pair of crossing sets $X, Y \in \calF$.
One can easily see that if $\calF$ is crossing,
$X, Y \in \calF$, $X \cap Y \ne \emptyset$, and $X \cup Y \ne V$,
then $X \cap Y, X \cup Y \in \calF$.
Finally, for a crossing family~$\calF$,
a function $f \colon \calF \to \R$ is called \emph{crossing submodular}
(on~$\calF$) if
$$
   f(X) + f(Y) \ge f(X \cap Y) + f(X \cup Y)
$$
holds for all $X, Y \in \calF$ such that $(X,Y)$ is a crossing pair.

We need some additional notation. For a set of nodes~$X$ denote the set of arcs entering
(resp. leaving) $X$ by $\deltain(X)$ (resp. $\deltaout(X)$).
Also $\gamma(X)$ (resp.~$\delta(X)$) will denote the set of arcs or edges having both
ends (resp. exactly one end) in $X$.

Put $\phi(X) := \abs{\deltain(X)}$. It is well-known
that $\phi$ is crossing submodular on $2^V$.
We consider the following subfamily of $2^V$:
$$
   \calF_1 := \setst{X \subseteq V}{\phi(X) = 1}.
$$
\begin{lemma}
   $\calF_1$ is a crossing family.
\end{lemma}
\begin{proof}
   Let $(X,Y)$ be a crossing pair of subsets of~$V$ such
   that $\phi(X) = \phi(Y) = 1$. Submodularity of $\phi$
   implies
   $$
      \phi(X \cap Y) + \phi(X \cup Y) \le \phi(X) + \phi(Y) = 2.
   $$
   On the other hand, since $X \cap Y \ne \emptyset$, $X \cup Y \ne V$
   and $G$ is strongly connected, one has $\phi(X \cap Y) \ge 1$ and
   $\phi(X \cup Y) \ge 1$. Therefore, $\phi(X \cap Y) = \phi(X \cup Y) = 1$
   and hence both $X \cap Y$ and $X \cup Y$ are members of $\calF_1$.
\end{proof}

\medskip

Consider a pair of multisets $S, T$ of nodes.
By an \emph{$S$--$T$ collection} we mean a collection of arc-disjoint paths in $G$
such that: (i)~each path of $\calP$ starts at a node in $S$ and ends at a node in $T$;
(ii)~for each $s \in S$ the number of paths from $\calP$ starting at~$s$
equals the multiplicity of $s$ in $S$;
(iii)~for each $t \in T$ the number of paths from $\calP$ ending at~$t$
equals the multiplicity of~$t$ in~$T$.

Let $x^+_i, y^+_i$ (resp.~$x^-_i, y^-_i$) be the ends of
$e^+_i$ (resp.~$e^-_i$). Consider the sets
$$
   V^+ := \setst{x^+_i, y^+_i}{1 \le i \le n}, \quad
   V^- := \setst{x^-_i, y^-_i}{1 \le i \le n}.
$$

Let $C$ be a cycle in $G + E$ that traverses each edge of~$E$.
Removing edges of~$E$ from~$C$ we split $C$ into a $V^+$--$V^-$ collection~$\calP_0$.
Consider an arbitrary index~$i$ and the sets $S := \set{x^+_1, y^+_1}$
and $T := \set{x^-_i, y^-_i}$. Suppose, there exists an $S$--$T$ collection.
Together with edges $e^+_1$ and $e^-_i$ these paths form a
cycle in $G + E$, contradicting the minimality of~$E$.
Therefore, no $S$--$T$ collection exists. Then, taking strong connectivity
of~$G$ into account, by a standard max-flow min-cut argument
there exists a set~$Z_i \in \calF_1$ such that $Z_i \cap S_1 = \emptyset$,
$T_i \subseteq Z_i$.

We start with sets $Z_1, \ldots, Z_n$ and
unite them to construct a collection of inclusion-wise maximum
sets $W_1, \ldots, W_m$. More precisely, let $H$ be an undirected
graph with nodes~$\set{1, \ldots, n}$. For each $1 \le i < j \le n$
we add an edge connecting nodes~$i$ and $j$ iff $Z_i \cap Z_j \ne \emptyset$.
Let $C_1, \ldots, C_m$ be the nodesets of connected components of~$H$.
For each $i$ put $W_i$ to be the union of $Z_j$, $j \in C_i$.
Clearly, $W_i \in \calF_1$ for all~$i$. From definition of $W_i$
it follows that $W_i$ are pairwise disjoint and $\bar W := \bigcup_i W_i$ covers
all nodes of~$V^-$.

For each~$i$ we say that nodes $x^+_i, y^+_i$ are \emph{mates}. In particular,
$x^+_i$ is the mate of~$y^+_i$, and $y^+_i$ is the mate of~$x^+_i$. Same terms are used
for $x^-_i$ and $y^-_i$. A simple inductive argument shows that for each $i$
and $t \in W_i$ the mate of~$t$ is also in $W_i$ and $S_1 \cap W_i = \emptyset$.

For a set $X \subseteq V$ put $n^+(X) := \abs{X \cap T^+}$
and $n^-(X) := \abs{X \cap T^-}$. It follows from the existence of~$\calP_0$ and max-flow
min-cut argument that
\begin{equation}
\label{eq:cut_sizes}
   n^+(X) \ge n^-(X) - 1 \qquad \mbox{for all $X \in \calF_1$}.
\end{equation}
In view of~\refeq{cut_sizes}, two cases are possible. First, one may have
\begin{equation}
\label{eq:no_small_cut}
   n^+(W_i) \ge n^-(W_i) \qquad \mbox{for all $1 \le i \le m$}.
\end{equation}
But since $W_i$ are disjoint,
\refeq{no_small_cut} implies that $n^+(\bar W) \ge n^-(\bar W)$. However,
all nodes in $T^-$ are covered by~$\bar W$
and at least two nodes in $T^+$ (namely, $x^+_1$ and $x^-_1$) are not covered
by~$\bar W$~--- a contradiction.

We may now assume that
\begin{equation}
\label{eq:small_cut}
   n^+(W_1) = n^-(W_1) - 1
\end{equation}
and $W_1$ covers the following pairs of mates in~$T^-$:
\begin{equation}
\label{eq:right_side}
   T_1 := \set{x^-_1, y^-_1, \ldots, x^-_q, y^-_q}.
\end{equation}

We claim that $q < n$. Suppose $q = n$, then $n^-(W_1) = n$.
However, $\set{x^+_1, y^+_1 } \cap W_1 = \emptyset$, thus $n^+(W_1) \le n - 2$.
This contradicts~\refeq{cut_sizes}.

Let $a_0$ denote the only arc in $G$ entering $W_1$ (recall that $W_1 \in \calF_1$).
$\calP_0$ contains a unique path ending in
each node of~$T_1$. Put $S_1 := V^+ \cap W_1$; by~\refeq{small_cut} $\abs{S_1} = \abs{T_1} - 1$
and there exists a unique node $v \in V^+ - S_1$ such that
$v$ is connected to some node in $T_1$, say~$x^-_1$, by the path~$P_0 \in \calP_0$
that crosses $\deltain(W_1)$ by~$a_0$.
We trace~$P_0$ starting from~$v$ until reaching~$a_0$;
let $R_0$ be the suffix of~$P_0$ starting with~$a_0$.

We construct a subcollection of $\calP_0$ as follows.
Initially, consider the node~$y^-_1$. It is connected by
the path~$P_1 \in \calP_0$ with the node in $S_1$ that we denote by~$x^+_1$.
If $y^+_1 \notin S_1$, then we stop.
Otherwise, $y^+_1$ is connected by the path~$Q_1 \in \calP_0$ with the node in $T_1$ that
we denote by $x^-_2$. We now consider its mate~$y^-_2$ and
proceed it the same way as we did for $y^-_1$.

In general, on the $i$-th step we consider the node~$y^-_i$
and find the corresponding path~$P_i \in \calP_0$. Let $x^+_i$
be the start node of~$P_i$. If $y^+_i \notin S_1$, we stop.
Otherwise, denote by~$Q_i \in \calP_0$ the path starting at~$y^+_i$.
Put $x^-_{i+1}$ to be the end node of~$Q_i$ and proceed with the next step.

This procedure eventually halts after, say, $l$ steps yielding a collection of paths
\begin{equation}
\label{eq:inner_paths}
   P_1, Q_1, \ldots, P_{l-1}, Q_{l-1}, P_l
\end{equation}
and a node $y^+_l \in T^+ - S_1$.
Note that all these paths are completely contained in $G[W_1]$.
Since $G$ is strongly connected, there exists a path~$Q_l$
from~$y^+_l$ to $x^-_1$. This path crosses $\deltain(W_1)$ and hence
$R_0$ is a suffix of~$Q_l$. Thus $Q_l$ is arc-disjoint from
all paths~\refeq{inner_paths}. Put
\begin{eqnarray*}
   \calP' & := & \set{P_1, Q_1, \ldots, P_{l-1}, Q_{l-1}, P_l, Q_l}, \\
   S' & := & \set{x^+_1, y^+_1, \ldots, x^+_l, y^+_l}, \\
   T' & := & \set{x^-_1, y^-_1, \ldots, x^-_l, y^-_l}
\end{eqnarray*}
Then $\calP'$ is an $S'$--$T'$ collection
that gives rise to a cycle in $G + E$ traversing some but not all edges of~$E$.
This, however, contradicts the minimality of~$E$. Proof of~\refth{two_edges} is now complete.

\section{Skew-Symmetric Graphs}
\label{sec:ss}

For bidirected graphs there is an alternative (and essentially
equivalent) language of \emph{skew-symmetric} graphs.
This section contains terminology and some basic facts
and explains the correspondence between skew-symmetric
and bidirected graphs. For a more detailed survey on skew-symmetric graphs,
see, e.g.,~\cite{tut-67,GK-96,GK-04,BK-05}.

A \emph{skew-symmetric graph} is a digraph~$G$ endowed with
two bijections $\sigma_V, \sigma_A$ such that: $\sigma_V$ is
an \emph{involution} on the nodes (i.e., $\sigma_V(v)\ne v$ and
$\sigma_V(\sigma_V(v)) = v$ for each node~$v$), $\sigma_A$ is an
involution on the arcs, and for each arc $a$ from $u$ to $v$,
$\sigma_A(a)$ is an arc from $\sigma_V(v)$ to $\sigma_V(u)$. For
brevity, we combine the mappings $\sigma_V, \sigma_A$ into one mapping
$\sigma$ on $V_G \cup A_G$ and call $\sigma$ the \emph{symmetry} (rather
than skew-symmetry) of $G$.
For a node (arc) $x$, its symmetric node (arc) $\sigma(x)$ is also
called the \emph{mate} of $x$, and we will often use notation with
primes for mates, denoting $\sigma(x)$ by $x'$.

Observe that if $G$ contains an arc $a$ from a node $v$
to its mate $v'$, then $a'$ is also an arc from $v$ to $v'$ (so the
number of arcs of $G$ from $v$ to $v'$ is even and these parallel
arcs are partitioned into pairs of mates).

The symmetry $\sigma$ is extended in a natural way to walks, paths, cycles,
and other objects in $G$. In particular, two walks are symmetric to
each other if the elements of the former are symmetric to those of the latter
and go in the reverse order:
for a walk $P = (v_0, a_1, v_1, \ldots, a_k, v_k)$, the symmetric
walk~$\sigma(P)$ is $(v'_k, a'_k, v'_{k-1}, \ldots, a'_1, v'_0)$.

Next we explain the correspondence between skew-symmetric and
bidirected graphs (cf.~\cite[Sec.~2]{GK-04}, \cite{BK-05}).
For sets $X, A, B$, we
use notation $X = A \sqcup B$ when $X = A \cup B$ and
$A \cap B= \emptyset$. Given a skew-symmetric
graph~$G$, choose an arbitrary partition $\pi=\set{V_1, V_2}$ of
$V_G$ such that $\sigma(V_1) = V_2$. Then $G$ and $\pi$ determine
the bidirected graph $\bar G$ with $V_{\bar G} := V_1$
whose edges correspond to the pairs of symmetric arcs in $G$. More precisely,
arc mates $a,a'$ of $G$ generate one edge $e$ of $\bar G$ connecting nodes
$u, v \in V_1$ such that: (i) $e$~goes from $u$ to $v$ if one of $a, a'$
goes from $u$ to $v$ (and the other goes from $v'$ to $u'$ in
$V_2$); (ii) $e$~leaves both $u,v$ if one of $a,a'$ goes from $u$
to $v'$ (and the other from $v$ to $u'$); (iii) $e$~enters both
$u, v$ if one of $a,a'$ goes from $u'$ to $v$ (and the other from
$v'$ to $u$). In particular, $e$ is a loop if $a, a'$ connect a pair
of symmetric nodes.

\begin{figure}[tb]
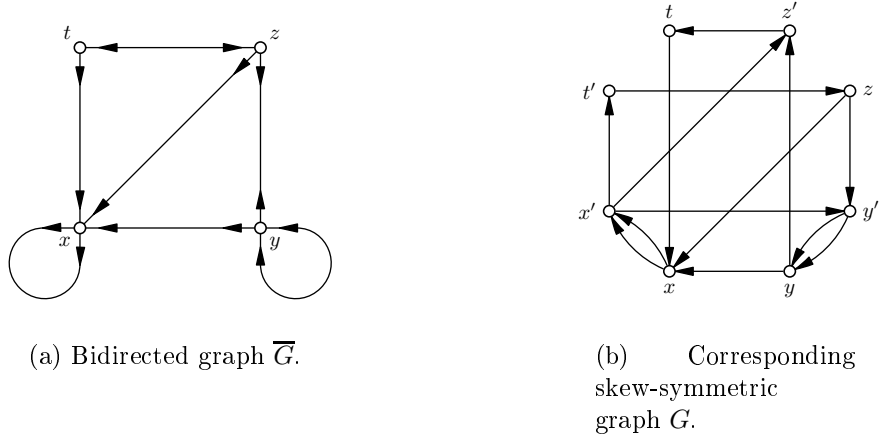

    \centering
    \subfigure[Bidirected graph~$\bar G$.]{
      \includegraphics{pics/examples.1}%
    }
    \hspace{3cm}%
    \subfigure[Corresponding skew-sym\-met\-ric graph~$G$.]{
      \includegraphics{pics/examples.2}%
    }
    \caption{Related bidirected and skew-symmetric graphs.}
    \label{fig:sk-bi}
\end{figure}

Conversely, a bidirected graph $\bar G$
determines a skew-sym\-met\-ric graph~$G$ with symmetry $\sigma$ as
follows. Take a copy $\sigma(v)$ of each element $v$ of $\bar V := V_{\bar G}$,
forming the set $\bar V' := \setst{\sigma(v)}{v\in \bar V}$. Now put
$V_G := \bar V \sqcup \bar V'$. For each edge $e$ of $\bar G$ connecting nodes
$u$ and $v$, assign two ``symmetric'' arcs $a,a'$ in $G$ so as to satisfy
(i)--(iii) above (where $u' = \sigma(u)$ and $v' = \sigma(v)$). An
example is depicted in Fig.~\ref{fig:sk-bi}.

Let $X$ be an arbitrary subset of nodes of a bidirected graph $\bar G$.
One can modify $\bar G$ as follows: for each node $v\in X$ and each edge
$e$ incident with $v$, reverse the direction of $e$ at $v$.
This transformation preserves the set of walks in $\bar G$ and
thus does not change the graph in essence.
We call two bidirected graphs $\bar G_1, \bar G_2$ \emph{equivalent}
if one can obtain $\bar G_2$ from $\bar G_1$ by applying a number
of described transformations.

\begin{remark}
   A bidirected graph generates one skew-symmetric graph,
   while a skew-symmetric graph generates a
   number of bidirected ones, depending on the partition $\pi$ of $V_G$.
   The latter bidirected graphs are equivalent.
\end{remark}

Also there is a correspondence between walks in $\bar G$ and walks in $G$.
More precisely, let $\tau$ be the natural mapping
of $V \cup A$ to $\bar V \cup \bar E$ (obtained by identifying the
pairs of symmetric nodes and arcs). Each walk
$P = (v_0, a_1, v_1, \ldots, a_k, v_k)$ in $G$ induces
the sequence
$$
   \tau(P) := (\tau(v_0), \tau(a_1), \tau(v_1), \ldots, \tau(a_k), \tau(v_k))
$$
of nodes and edges in $\bar G$. One can easily check that $\tau(P)$ is a walk in $\bar G$
and $\tau(P') = \tau(P)^R$ (where $W^R$ stands for the bidirected walk
obtained by passing $W$ in opposite direction).
Moreover, for any walk $\bar P$ in $\bar G$
there is exactly one pre-image $\tau^{-1}(\bar P)$.

\section{Regular Reachability and Barriers}
\label{sec:barrier}

A path in a skew-symmetric graph is called \emph{regular} if it does
not contain a pair of symmetric arcs (while symmetric nodes are allowed).
This notion plays an important role since regular paths in a skew-symmetric
graph~$G$ are exactly the images of paths in the corresponding bidirected
graph~$\bar G$. In this section we
state a criterion for the existence of a regular path connecting a pair
of symmetric nodes in a skew-symmetric graph.

Consider a skew-symmetric graph~$G$.
Let $\tau = (V_\tau, a_\tau)$, $V_\tau \subseteq V_G$, $a_\tau \in A_G$ be a pair such that:
(i) $V_\tau' = V_\tau$; (ii) $a_\tau \in \deltain(V_\tau)$;
(iii) every node in $V_\tau$ is reachable from the head of $a_\tau$ by a regular path
in~$G[V_\tau]$. Then we $\tau$ is called a \emph{bud}.

\begin{figure}[tb]
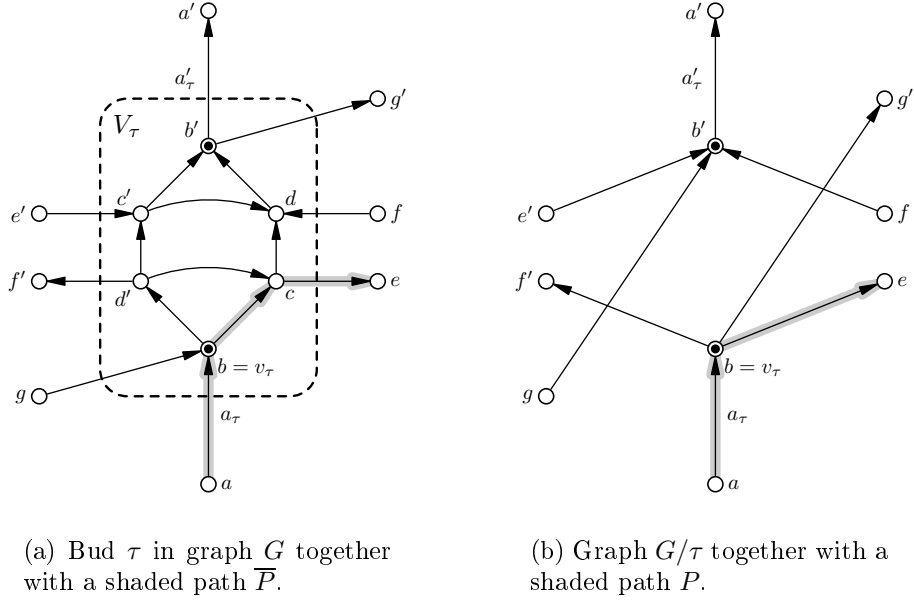

   \centering
   \subfigure[Bud~$\tau$ in graph~$G$ together with a shaded path~$\bar P$.]{
      \includegraphics{pics/trimming.1}%
      \label{fig:bud}
   }
   \hspace{1cm}%
   \subfigure[Graph~$G / \tau$ together with a shaded path $P$.]{
      \includegraphics{pics/trimming.2}%
   }
   \caption{Buds, trimming, and path restoration.
   Base and antibase nodes $b, b'$ are marked. Path $\bar P$ is a preimage of~$P$.}
   \label{fig:trimming}
\end{figure}

Let $v_\tau$ denote the head node of~$a_\tau$.
The arc~$a_\tau$ (resp. node~$v_\tau$)
is called the \emph{base arc} (resp. \emph{base node}) of $\tau$,
arc~$a_\tau'$ (resp. node~$v_\tau'$) is called the \emph{antibase arc} (resp.
\emph{the antibase node}) of $\tau$. For an arbitrary bud~$\tau$ we denote its
set of nodes by $V_\tau$, base arc by $a_\tau$, and base node by $v_\tau$.
An example of a bud is given in \reffig{bud}.

Consider an arbitrary bud~$\tau$ in a skew-symmetric graph~$G$.
By \emph{trimming~$\tau$} we mean the following transformation of $G$:
(i) all nodes in $V_\tau - \set{v_\tau, v_\tau'}$ and arcs in $\gamma(V_\tau)$
are removed;
(ii) all arcs in $\deltain(V_\tau) - \set{a_\tau}$ are transformed into
arcs entering $v_\tau'$ (the tails of these arcs are not changed);
(iii) all arcs in $\deltaout(V_\tau) - \set{a_\tau'}$ are transformed into
arcs leaving $v_\tau$ (the heads of these arcs are not changed).
The resulting skew-symmetric graph is denoted by~$G / \tau$.
Thus, each arc of the original graph $G$ not belonging to $\gamma(V_\tau)$
has its \emph{image} in the trimmed graph~$G / \tau$.
\reffig{trimming} gives an example of bud trimming.

Let $P$ be a regular path in $G / \tau$. One can lift this path to $G$ as follows:
if $P$ does not contain neither $a_\tau$, nor $a_\tau'$ leave~$P$ as it is.
Otherwise, consider the case when $P$ contains~$a_\tau$ (the symmetric case is analogous).
Split $P$ into two parts: the part $P_1$ from the beginning of $P$ to $v_\tau$ and
the part~$P_2$ from $v_\tau$ to the end of $P$. Let $a$ be the first arc of $P_2$. The arc~$a$
leaves $v_\tau$ in $G / \tau$ and thus corresponds to some arc $\bar a$ leaving $V_\tau$ in $G$
($\bar a \ne a_\tau'$). Let $u \in V_\tau$ be the tail of $a$ in $G$ and $Q$ be
a regular $v_\tau$--$u$ path in $G[V_\tau]$ (existence of $Q$ follows from definition of bud).
Consider the path $\bar P := P_1 \circ Q \circ P_2$ (here $U \circ V$ denotes the
path obtained by concatenating $U$ and $V$). One can easily show that $\bar P$ is regular.
We call $\bar P$ a \emph{preimage of~$P$} (under trimming $G$ by $\tau$).
Clearly, $\bar P$ is not unique.
An example of such path restoration is shown in \reffig{trimming}: the shaded path~$\bar P$
on the left picture corresponds to the shaded path~$P$ on the right picture.

\begin{figure}[tb]
   \centering
   \includegraphics{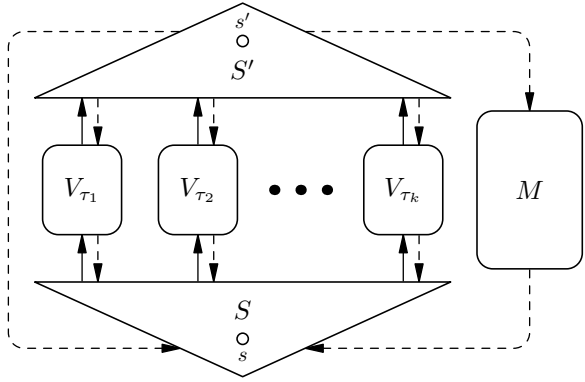}%
   \caption{
      An $s$-barrier. Solid arcs should occur exactly once,
      dashed arcs may occur arbitrary number of times (including zero).
   }
   \label{fig:barrier}
\end{figure}

Let $G$ be a skew-symmetric graph with a designated node $s$.
Suppose we are given a collection of buds $\tau_1, \ldots, \tau_k$ in $G$ together with
node sets $S$ and $M$. Additionally, suppose the following properties hold:
(i) collection $\set{S, S', M, V_{\tau_1}, \ldots, V_{\tau_k}}$ forms a partition of~$V_G$
with $s \in S$;
(ii) no arc goes from $S$ to $S' \cup M$;
(iii) no arc connects distinct sets $V_{\tau_i}$ and $V_{\tau_j}$;
(iv) no arc connects $V_{\tau_i}$ and $M$;
(v) the arc~$a_{\tau_i}$ is the only one going from $S$ to $V_{\tau_i}$.
Then we call the tuple $\calB = (S, M; \tau_1, \ldots, \tau_k)$
an $s$-\emph{barrier} (\cite{GK-96}, see~\reffig{barrier}
for an example).

\begin{theorem}[Barrier Theorem, \cite{GK-96}]
\label{th:regular_reach_criterion}
   There exists a regular $s$--$s'$ path in a skew-symmetric graph~$G$
   iff there is no $s$-barrier in $G$.
\end{theorem}

\section{Proof of \refth{bidir_decomp}}
\label{sec:bidir_decomp_proof}

By an inductive argument it is sufficient to prove that given a
strongly connected bidirected graph~$\bar G$ and its strongly connected proper subgraph~$\bar H$
one can extend $\bar H$ to a strongly connected graph by adding one or two edge-disjoint
ears of~$\bar H$ w.r.t.~$\bar G$. Moreover, one may assume that no single-ear step is possible
at the moment and prove that a double-ear step can be performed in this case.

Consider skew-symmetric graphs~$G$ and $H$ that are
related to $\bar G$ and $\bar H$ respectively. Let $a_0$ be an arc from $A_G - A_H$
that has its tail node~$u_0$ in~$V_H$ (such arc exists due to
connectivity of underlying undirected graphs of~$\bar G$ and $\bar H$).
Since $\bar G$ is strongly connected there exists a regular
cycle~$C_0$ in~$G$ passing through~$a_0$.
We follow along this cycle
starting from~$a_0$ until reaching the nodeset of~$H$. This way, we construct
a path~$P_0$ in~$G$ from~$u_0 \in V_H$ to, say, $v_0 \in V_H$.
The image of~$P_0$ in $\bar H$ forms an ear w.r.t.~$\bar G$.

By assumption that no single-ear step is currently possible, one has no regular path in~$H$
from $v_0$ to $u_0$. To apply~\refth{regular_reach_criterion} we construct an auxiliary skew-symmetric
graph~$H_0$ from~$H$ by adding a pair of symmetric nodes $s,s'$ and arcs~$(s,v_0)$, $(s,u_0')$,
$(v_0',s')$, $(u_0,s')$. It follows that no regular $s$--$s'$ path exists in~$H_0$ and thus
there exists an $s$-barrier $\calB_0 = (\set{s} \cup A, M; \tau_1, \ldots, \tau_k)$ in~$H_0$
where $A, M \subseteq V_H$ and $\tau_i$ are buds in $H_0$.

\begin{lemma}
   $\calB := (A, \emptyset; \tau_1, \ldots, \tau_k)$ is a $v_0$-barrier in~$H$.
\end{lemma}
\begin{proof}
   First, suppose that $\tau_i$ is not a bud in $H$. This is only possible if
   the tail of its base arc~$a_{\tau_i}$ is~$s$. Hence,
   \begin{equation}
   \label{eq:not_connected}
      \deltain_H(V_{\tau_i}) = \deltaout_H(V_{\tau_i}) = \emptyset,
   \end{equation}
   that a contradiction with connectivity of the underlying undirected graph of~$\bar H$.
   Therefore, all $\tau_i$ are also buds in $H$. To see that $v_0 \in A$ note
   that the only other possibility for $v_0$ is to be the base node of some bud~$\tau_i$.
   This, however, would again imply~\refeq{not_connected} and hence is not possible.
   We also prove that $M = \emptyset$. Indeed, if $\deltain(M) = \deltaout(M) = \emptyset$,
   then the underlying undirected graph of~$\bar H$ is not connected. In case there exists
   an arc leaving $M$, from definition of barrier it follows
   that no regular cycle in~$H$ can pass through this arc~--- again a contradiction.
\end{proof}

Consider the graph $H_1 := H / \tau_1 / \ldots / \tau_k$ obtained from~$H$ by trimming
all buds of~$\calB$. Put $Z := A \cup \set{v_{\tau_1}, \ldots, v_{\tau_k}}$ and
consider the bidirected graph~$\bar H_1$ corresponding to~$H_1$ under partition $\set{Z, Z'}$
of~$V_{H_1}$ (see~\refsec{intro}).
Since no arc in~$H_1$ connects the sets~$Z$ and $Z'$, all edges of~$\bar H_1$
are standard, so we may regard $\bar H_1$ as a digraph
isomorphic to $H_1[Z]$. As long as this is not confusing, we make no distinction between
$\bar H_1$ and $H_1[Z]$.

\begin{lemma}
\label{lm:trimmed_conn}
   $H_1[Z]$ is strongly connected.
\end{lemma}
\begin{proof}
   The connectivity of the underlying undirected graph follows from this property of~$\bar H$.
   Consider an arbitrary arc~$a$ of~$H_1[Z]$.
   Consider a regular cycle~$C$ passing through~$a$ in $H$; $C$ remains a regular
   cycle under trimming of all buds in~$\calB$. The image of~$C$ under these trimmings gives rise
   to a cycle in~$H_1[Z]$ that passes through~$a$, as required.
\end{proof}

\begin{figure}[bt]
    \centering
    \includegraphics{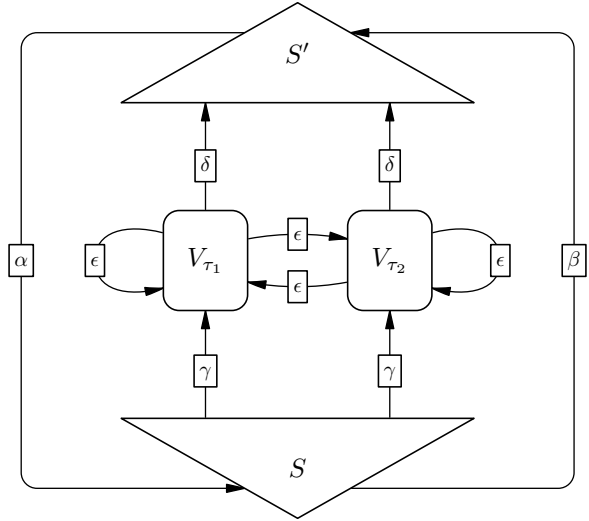}%
    \caption{Possible types of ears.}
    \label{fig:ears}
\end{figure}

Recall that we originally had the arc~$a_0 \in A_G - A_H$ and the regular cycle~$C_0$
passing through~$a_0$.
We drop all arcs of~$C_0$ that belong to $A_H$ and thus split~$C_0$ into
a collection of ears of~$H$ w.r.t.~$G$. Consider an arbitrary such
ear~$P$; let $u$ be its start node, and $v$ be its end node. We call~$P$
\begin{enumerate}
   \item \emph{$\alpha$-ear} if $u \in S'$, $v \in S$;
   \item \emph{$\beta$-ear} if $u \in S$, $v \in S'$;
   \item \emph{$\gamma$-ear} if $u \in S$, $v \in V_{\tau_i}$ for some~$i$;
   \item \emph{$\delta$-ear} if $u \in V_{\tau}$ for some~$i$, $v \in S'$;
   \item \emph{$\epsilon$-ear} if $u \in V_{\tau_i}$, $v \in V_{\tau_j}$ for some $i, j$ (possibly $i = j$).
\end{enumerate}
These five cases are depicted in \reffig{ears}.

\begin{lemma}
   Each ear obtained from~$C_0$ belongs to one these five classes.
\end{lemma}
\begin{proof}
   Let $P$ be an ear not falling into one of these classes.
   Due to symmetry, it is sufficient to consider the following two cases:
   (i)~$u, v \in S$;
   (ii)~$u \in V_{\tau_i}$ for some~$i$ and $v \in S$.
   We argue that $u$ is reachable from $v$ by a regular path in~$H$
   and hence the image of $P$ in $\bar G$ is an ear that can be added to $\bar H$ without
   loss of strong connectivity. This contradicts the assumption that
   no single-ear step is currently possible.

   Indeed, in (i) \reflm{trimmed_conn} implies that $u$ is reachable from $v$ by a
   regular path in $H_1$. By a standard restoration procedure this path can be extended
   to a regular $v$--$u$ path in $H$.
   In (ii) $v_{\tau_i}$ is reachable from $v$ by a regular path in $H_1$.
   Applying restoration procedure and adding a regular
   $v_{\tau_i}$--$u$ path in $H[V_{\tau_i}]$ one again gets a regular
   $v$--$u$ path in~$H$.
\end{proof}

Next, we consider the sequence of ears
\begin{equation}
\label{eq:ears_of_cycle}
   P_0, P_1, \ldots, P_m
\end{equation}
obtained from $C_0$ and construct a collection of nonstandard edges~$E$ such
that $\bar H_1 + E$ is strongly connected.
First consider the set of $\alpha$-ears in~\refeq{ears_of_cycle}.
Each such ear~$P$ (in particular, $P_0$) goes from a node $u \in S'$ to a node $v \in S$.
Construct an edge $\set{u, v}$ (called \emph{backward})
that enters both of its ends and assign the ear~$P$ to this edge.

The subsequence of $\alpha$-ears splits \refeq{ears_of_cycle} into maximal parts
without $\alpha$-ears.
Let $P_i, \ldots, P_j$ be any of these parts. The part gives rise
to a nonstandard edge as follows.
In case~$P_i$ is $\beta$- or $\gamma$-ear
put $x$ to be the start node of~$P_i$. Otherwise ($P_i$ is $\delta$- or $\epsilon$-ear)
put $x$ to be the base node of the bud containing the start node of~$P_i$.
Similarly, consider $P_j$.
In case~$P_j$ is $\beta$- or $\delta$-ear
put $y$ to be the end node of~$P_i$. Otherwise ($P_j$ is $\gamma$- or $\epsilon$-ear)
put $y$ to be the antibase node of the bud containing the end node of~$P_j$.
Construct an edge $\set{x, y}$ (called \emph{forward})
that leaves both of its ends and assign the sequence
of ears $P_i, \ldots, P_j$ to this edge.

As a result, we get a collection of nonstandard edges~$E$.
All these edges belong to a cycle in~$\bar H_1 + E$ obtained from~$C_0$ as follows:
\begin{numitem}
\label{eq:cycle_image}
   All arcs in $\gamma(V_{\tau_i})$, $i = 1, \ldots, k$ are dropped.
   Each $\alpha$-ear~$Q$ in \refeq{ears_of_cycle} is replaced by the arc
   corresponding to the backward edge assigned to~$Q$. Each maximal sequence $P_i, \ldots, P_j$
   of $\beta$-, $\gamma$-, $\delta$-, and $\epsilon$-ears is replaced
   by the arc corresponding to the forward edge assigned to $P_i, \ldots, P_j$.
   Finally, the bidirected image in $\bar H_1 + E$ is taken by merging mates of nodes and arcs.
\end{numitem}

Hence, $\bar H_1 + E$ is strongly connected.
Now \refth{two_edges} implies the existence of a pair of edges $e_1, e_2 \in E$
(where $e_1$ is forward and $e_2$ is backward) such that $\bar H_1 + e_1 + e_2$ is strongly connected.
Our final task is to replace these edges by a pair of ears of~$\bar H$ w.r.t.~$\bar G$.

A trivial part is to deal with~$e_2$ since it corresponds
to a single ear in~\refeq{ears_of_cycle}.
In contrast, $e_1$ may correspond to a number of ears.
We first prove the following auxiliary statement:

\begin{lemma}
\label{lm:five_nodes}
   Consider an arbitrary strongly connected digraph and
   five nodes $a, b, x, y, z$ in it. Suppose there exists
   an $\set{a,b}$--$\set{x,y}$ collection of paths.
   Then, there exists an $\set{a,b}$--$\set{x,z}$ or
   $\set{a,b}$--$\set{z,y}$ collection.
\end{lemma}
\begin{proof}
   We may assume that there exist an $a$--$x$ path~$P$
   and a $b$--$y$ path~$Q$ that are arc-disjoint.
   Consider an arbitrary $a$--$z$ path~$R$. We follow
   it backwards starting from~$z$ and stop either when reaching~$a$ or encountering
   an arc from~$A_P$ or $A_Q$. If $a$ is reached,
   then $\set{Q,R}$ is a desired $\set{a,b}$--$\set{z,y}$
   collection. If an arc from~$P$ is encountered, then we get an
   $\set{a,b}$--$\set{z,y}$ collection by taking
   path~$Q$ and parts of paths~$P, R$. Finally, if
   an arc from~$P$ is encountered, an $\set{a,b}$--$\set{x,z}$ collection is obtained
   by taking path~$P$ and parts of paths~$Q, R$.
\end{proof}

\medskip

To complete the proof we now proceed iteratively as follows. We maintain a pair
of nonstandard edges $e_1, e_2$ ($e_1$ is forward, $e_2$ is backward).
Edge~$e_2$ is assigned a $\alpha$-ear from~\refeq{ears_of_cycle};
let us denote this ear by~$Q$. Edge~$e_1$ is assigned a sequence $P_i, \ldots, P_j$ of
ears from~\refeq{ears_of_cycle}.
The following invariant holds: there exists a regular cycle~$C$
in the skew-symmetric graph $H + (Q + Q') + (P_i + P_i') + \ldots + (P_j + P_j')$
that passes through all arcs of $Q, P_i, \ldots, P_j$ in this order.
Moreover, $C$ gives rise to a cycle~$\bar C_1$ in $\bar H_1 + e_1 + e_2$ according
to~\refeq{cycle_image}.

Let $\set{a,b}$ be the multiset of ends
of~$e_2$ and $\set{x,y}$ be the multiset of ends
of $e_1$. Due to symmetry, we may assume that $x$ is the start node of~$P_i$.
(Hereinafter we identify nodes
of $\bar H_1$ with those of $H_1[Z]$ and $H[Z]$.)
By dropping edges $e_1, e_2$ from~$\bar C_1$ one gets
an $\set{a,b}$--$\set{x,y}$ collection of paths in $H_1$.
In case $i = j$, a unique ear corresponds to $e_1$ and hence we are done.
Otherwise we change $e_1, i, j$ so as to reduce the number of ears assigned to~$e_1$.
Consider $P_i$; it cannot be a $\beta$- or $\delta$-ear since that would imply $i = j$.
Hence, two cases are possible.

If $P_i$ is a $\gamma$-ear then put $z$ to be base node of the bud
containing the end node of~$P_i$. Apply~\reflm{five_nodes} and replace
$\set{x,y}$ by either $\set{z,y}$ or $\set{x,z}$.
In the former case put $i := i + 1$, in the latter put $j := i$.
Also, update the cycle~$C$ and the edge~$e_1$ to reflect the changes in its ends
and proceed with the next iteration.

Now suppose $P_i$ is $\epsilon$-ear with the start node in the nodeset of
a certain bud, say~$\tau$. In this case $x = v_\tau$.
The cycle~$C$ enters $V_\tau$ by the arc~$a_\tau$,
uses some arcs from $\gamma_H(V_\tau)$, and then leaves~$V_\tau$
by~$P_i$. We make sure that $P_i$ is the only ear assigned to $e_1$ that leaves $V_\tau$.
If it is not true then we replace $i$ by the largest index~$k$ in the range $i, \ldots, j$
such that $P_k$ leaves $V_\tau$. The cycle~$C$ and the edge~$e_1$ are updated accordingly.

Then, let $\eta$ be the bud whose nodeset contains the end node of~$P_i$;
put $z := v_{\eta}$. Like earlier, we apply~\reflm{five_nodes} and
and replace $\set{x,y}$ by either $\set{z,y}$ or $\set{x,z}$.
In the former case put $i := i + 1$, in the latter put $j := i$.
As before, update the cycle~$C$ and the edge~$e_1$ to reflect the changes in its ends $x, y$
and proceed with the next iteration.

Once iterations are complete, we get a  single ear~$P_i$ is assigned to $e_1$.
The bidirected images of $P_i$, $Q$ in $\bar G$ form the desired pair of
ears of~$\bar H$ w.r.t.~$\bar G$. The proof of~\refth{bidir_decomp} is now complete.

\section{Application to Matching Covered Graphs}
\label{sec:applications}

Recall \cite{LP-86}
that a \emph{perfect matching}~$M$ in an undirected graph~$G$ is a set of edges
such that each node $v \in V_G$ is incident with exactly one edge in~$M$.
An undirected graph is called \emph{matching covered} if every edge $e \in E_G$
is contained in a perfect matching. An  path in $G$ is called
\emph{alternating} w.r.t.~$M$ if it consists of an alternating sequence
of edges from $M$ and $E_G - M$.

A subgraph $H$ of $G$ is called \emph{elastic} (w.r.t.~$G$) if $G[V_G - V_H]$ has
a perfect matching. By an \emph{ear} of $H$ w.r.t.~$G$ we mean
a simple path~$P$ of odd length in $G$ such that:
(i)~ends of~$P$ are distinct and are contained in~$V_H$;
(ii)~no inner node of~$P$ is contained in~$V_H$;
(iii)~ $E_P \cap E_H = \emptyset$.
The result of adding $P$ to $H$ is denoted by $H + P$ and is
defined in a natural way.

An \emph{ear decomposition} of a matching covered graph $G$ starting from its
elastic matching covered subgraph~$H$ is a sequence of elastic (w.r.t.~$G$)
matching covered subgraphs of~$G$
$$
   H = G_0, G_1, \ldots, G_{k-1}, G_k = G,
$$
where $G_{i+1}$ is obtained
from $G_i$ by adding a single ear of~$G_i$ w.r.t.~$G$
or a node-disjoint pair of such ears ($0 \le i < k$).

We use~\refth{bidir_decomp} to derive the following result of Lov\'asz and Plummer:
\begin{theorem}
   For any matching covered graph~$G$ and an arbitrary elastic subgraph~$H$ of $G$
   there exists an ear decomposition of~$G$ starting from~$H$.
\end{theorem}
\begin{proof}
   It is sufficient to prove that for a matching covered graph $G$ and its
   elastic matching covered proper subgraph $H$ the latter one can be extended
   to an elastic matching covered graph by adding one or two ears of~$H$ w.r.t.~$G$.

   Consider a perfect matching~$M_G$ in~$G$ such that $M_H := M_H \cap E_H$
   is a perfect matching in~$H$ (existence of~$M_G$ follows from elasticity of~$H$).
   Then $G$, $M_G$ generate
   the bidirected graph~$\bar G$ as follows. Each edge $e \in E_G - M_G$
   is directed so as to leave both of its ends. Each edge $e = \set{u,v} \in M$
   is transformed into a pair of parallel edges $e_1, e_2$ connecting nodes $u, v$.
   The former one enters $u$ and $v$; the latter one leaves $u$ and $v$.
   Edges~$e_2$ are called \emph{auxiliary}. A similar construction applied to $H, M_H$
   yields the bidirected subgraph~$\bar H$ of~$\bar G$.

   We prove that $\bar G$ is strongly connected (the same argument also applies to $\bar H$).
   For each $e \in E_G$ the edges $e_1, e_2$ form a cycle in~$\bar G$. So it remains
   to consider edges $e \in E_G - M_G$. From definition of matching covered
   graph and simple facts regarding perfect matchings it follows that there
   exists an alternating cycle in $G$ w.r.t.~$M_G$ that passes through~$e$.
   This cycle in $G$ gives rise to a desired cycle in $\bar G$ passing through~$e$.

   Consider an arbitrary ear~$\bar P$ of $\bar H$ w.r.t.~$\bar G$
   and its image~$P$ in~$G$ (obtained by dropping directions of edges
   and merging $e_1, e_2$ into $e$, where $e \in M_G$). Suppose~$\bar P$ contains
   an auxiliary edge~$e_2$ (corresponding to the edge~$e = \set{u,v} \in M_G$).
   It follows that both ends of~$\bar P$ are
   contained in the set $\set{u,v}$. Hence, $\set{u,v} \subseteq V_H$ and $e \in M_H$.
   Thus, $e_2 \in E_{\bar H}$, which is a contradiction.
   It is now easy to see that~$P$ is an alternating path in~$G$ w.r.t.~$M_G$
   with first and last edges in $E_G - M_G$. Therefore, it has an odd length.

   We apply~\refth{bidir_decomp} to $\bar H, \bar G$ to get a collection~$\bar \calP$
   of at most two edge-disjoint ears of~$\bar H$ w.r.t.~$\bar G$ such that
   adding all ears of~$\bar \calP$ to $\bar H$ one gets a strongly connected graph
   $\bar H' := \bar H + \bar \calP$. We may assume that $\abs{\bar \calP}$ is minimal and hence
   there exists a cycle $\bar C$ in
   $\bar H'$ that passes through all ears from~$\bar \calP$. The image~$C$
   of $\bar C$ in $G$ is an alternating cycle w.r.t.~$M_G$. Each alternating
   cycle is simple and thus all nodes of ears in~$\bar \calP$ are distinct, as
   required.

   It remains to show that the graph $H'$ (obtained from~$H$ by adding the images of
   ears from~$\bar \calP$) is elastic and matching covered. The former
   property follows from the fact that $M_G \cap \gamma(V_G - V_{H'})$ is a perfect matching
   in $G[V_G - V_{H'}]$. The latter property is due to the strong connectivity of~$\bar H'$.
\end{proof}

\section{Acknowledgments}
The author is thankful to Alexander Karzanov for constant attention,
collaboration, and many insightful discussions.

\nocite{*}
\bibliographystyle{plain}
\bibliography{main}

\end{document}

%% file: abstract.tex
Bidirected graphs (earlier studied by Edmonds, Johnson and,
in equivalent terms of skew-symmetric graphs, by
Tutte, Goldberg, Karzanov, and others) proved to be a
useful unifying language for describing both
flow and matching problems.  In this paper we extend the notion
of ear decomposition to the class of strongly connected bidirected graphs.
In particular, our results imply Two Ear Theorem on matching
covered graphs of Lov\'asz and Plummer.
The proofs given here are self-contained except for standard
Barrier Theorem on skew-symmetric graphs.